\documentclass[10pt]{article}
\usepackage{cite}
\usepackage{mathrsfs}
\usepackage{amsfonts}
\usepackage{amsmath}
\usepackage{amsfonts,amssymb,color}
\usepackage{dsfont}
\usepackage{curves}
\usepackage{mathrsfs}
\usepackage{pifont}
\usepackage{amssymb}
\allowdisplaybreaks

\numberwithin{equation}{section}

\date{}

\def\BigRoman{\uppercase\expandafter{\romannumeral\number\count 255 }}
\def\Romannumeral{\afterassignment\BigRoman\count255=}

\begin{document}
\title{Spectral radii and star-factors with large components
}
\author{\small Zhiren Sun$^{1,2}$, Sizhong Zhou$^{3}$\footnote{Corresponding
author. E-mail address: zsz\_cumt@163.com (S. Zhou)}\\
\small $1$. Basic Education Department, Xinjiang University of Political Science and Law,\\
\small Tumushuke, Xinjiang 843900, China\\
\small $2$. School of Mathematical Sciences, Nanjing Normal University,\\
\small Nanjing, Jiangsu 210023, China\\
\small $3$. School of Science, Jiangsu University of Science and Technology,\\
\small Zhenjiang, Jiangsu 212100, China\\
}

\maketitle
\begin{abstract}
\noindent Let $G$ be a connected graph with $n$ vertices. The isolated toughness of $G$, denoted by $I(G)$, is defined by $I(G)=\min\left\{\frac{|S|}{i(G-S)}:S\subseteq V(G) \ \mbox{and} \ i(G-S)\geq2\right\}$
if $G$ is not complete, or $I(G)=+\infty$ if $G$ is complete. A graph $G$ is called isolated $r$-tough if $I(G)\geq r$. A spanning subgraph $H$ of $G$ is called a $\{K_{1,j}:m\leq j\leq2m\}$-factor of $G$ if
every component of $H$ is isomorphic to an element of $\{K_{1,j}:m\leq j\leq2m\}$. Let $\rho(G)$, $q(G)$ and $\mu(G)$ denote the adjacency spectral radius, the signless Laplacian spectral radius and the distance
spectral radius of $G$, respectively. Let $m$ and $b$ be two positive integers with $m\geq2$. In this paper, we first establish a lower bounds on the adjacency spectral radius of a connected isolated
$\frac{mb-1}{b}$-tough graph $G$ to guarantees that $G$ contains a $\{K_{1,j}:m\leq j\leq2m\}$-factor. Second, we establish a lower bounds on the signless Laplacian spectral radius of a connected isolated
$\frac{mb-1}{b}$-tough graph $G$ to ensures that $G$ contains a $\{K_{1,j}:m\leq j\leq2m\}$-factor. Finally, we create an upper bounds on the distance spectral radius of a connected isolated $\frac{mb-1}{b}$-tough
graph $G$ with a $\{K_{1,j}:m\leq j\leq2m\}$-factor. Furthermore, we construct some extremal graphs to claim that all the bounds obtained in this paper are sharp.
\\
\begin{flushleft}

(2020) Mathematics Subject Classification: 05C50, 05C70
\end{flushleft}
\end{abstract}

\medskip

\section{Introduction}

Throughout this paper, we deal only with finite, undirected and simple graphs. Let $G$ be a graph with vertex set $V(G)$ and edge set $E(G)$. The order of $G$, denoted by $n$, is defined as the number of its vertices. That
is to say, $|V(G)|=n$. The degree of a vertex $v$ in $G$, denoted by $d_G(v)$, is defined as the number of vertices adjacent to $v$ in $G$. A vertex $v$ with $d_G(v)=0$ in $G$ is called an isolated vertex of $G$. The number
of isolated vertices in $G$ is denoted by $i(G)$. For a subset $S\subseteq V(G)$, the subgraph of $G$ induced by $S$ is denoted by $G[S]$, and we write $G-S=G[V(G)\setminus S]$. The cycle, the path, the star and the complete
graph of order $n$ are denoted by $C_n$, $P_n$, $K_{1,n-1}$ and $K_n$, respectively. Given two vertex-disjoint graphs $G_1$ and $G_2$, their union is denoted by $G_1\cup G_2$. For any positive integer $t$, we denote by $tG$
the disjoint union of $t$ copies of $G$. The join $G_1\vee G_2$ is obtained from $G_1\cup G_2$ by joining every vertex of $G_1$ with every vertex of $G_2$ by an edge.

Yang, Ma and Liu \cite{YML} first introduced the concept of isolated tough, which is denoted by $I(G)$. The isolated tough of $G$ is defined as
$$
I(G)=\min\left\{\frac{|S|}{i(G-S)}:S\subseteq V(G) \ \mbox{and} \ i(G-S)\geq2\right\}
$$
if $G$ is not complete, or $I(G)=+\infty$ if $G$ is complete. A graph $G$ is called isolated $t$-tough if its isolated tough $I(G)\geq t$. A graph $G$ is isolated $t$-tough if and only if $t\cdot i(G-S)\leq|S|$ for every subset
$S\subseteq V(G)$ with $i(G-S)\geq2$.

Given a graph $G$ with $V(G)=\{v_1,v_2,\ldots,v_n\}$, the distance between two vertices $v_i$ and $v_j$, denoted by $d_{ij}$, is defined as the minimum of the length of the $v_i-v_j$ paths in $G$. The Wiener index, denoted by
$W(G)$, of a connected graph $G$ with $n$ vertices is defined as the sum of all distances in $G$. That is to say, $W(G)=\sum\limits_{i<j}d_{ij}$. Let $A(G)$ denote the adjacency matrix of $G$. The largest eigenvalue of $A(G)$
is called the adjacency spectral radius of $G$, which is denoted by $\rho(G)$. The signless Laplacian matrix of $G$ is defined by $Q(G)=D(G)+A(G)$, where $D(G)=\mbox{diag}(d_G(v_1),d_G(v_2),\ldots,d_G(v_n))$ is the degree
diagonal matrix of $G$. The largest eigenvalue of $Q(G)$ is called the signless Laplacian spectral radius of $G$, which is denoted by $q(G)$. The distance matrix of $G$, denoted by $\mathcal{D}(G)$, is defined by
$\mathcal{D}(G)=(d_{ij})_{n\times n}$. The largest eigenvalue of $\mathcal{D}(G)$ is called the distance spectral radius of $G$, which is denoted by $\mu(G)$.

A spanning subgraph of a graph is called a factor of the graph. Let $a$ and $b$ be two positive integers with $a\leq b$. Then a spanning subgraph $F$ with $a\leq d_F(v)\leq b$ of $G$ is called an $[a,b]$-factor of $G$. In
particular, a $[1,1]$-factor is simply called a 1-factor, and an $[a,b]$-factor is a $[1,k]$-factor if $a=1$ and $b=k$. Note that all of these notions are relate to the degree of vertices. Hence, they are often regarded as
``degree factors". Lots of results on this topic can be referred to \cite{AK,AGS,GWC,Zr,Wa}.

On the other hand, when focusing on components of a factor, we can lead to the concept of ``component factors". A set of connected graphs is denoted by $\mathcal{H}$. Then a spanning subgraph $H$ of $G$ is called an
$\mathcal{H}$-factor of $G$ if every component of $H$ is isomorphic to an element of $\mathcal{H}$. For an integer $m$ with $m\geq2$, an $\mathcal{H}$-factor is called a $\{K_{1,j}:1\leq j\leq m\}$-factor, which is also called
a star-factor. For $\mathcal{H}=\{P_m,P_{m+1},\ldots\}$, an $\mathcal{H}$-factor is called a $P_{\geq m}$-factor. In particular, a $\{K_{1,1}\}$-factor (or $\{P_2\}$-factor) is also called a 1-factor, a $P_{\geq m}$-factor is
a $[1,2]$-factor, and a $\{K_{1,j}:1\leq j\leq m\}$-factor is a $[1,m]$-factor.

Dai \cite{Dai}, Liu and Pan \cite{LP}, Zhang \cite{Zb} proposed some sufficient conditions for graphs with path-factors. Egawa and Furuya \cite{EF} investigated the existence of $K_{1,1}$-factors in star-free graphs. Tutte
\cite{Tt} verified a result, which provides the existence of a $\{K_{1,1},C_m:m\geq3\}$-factor based on the number of isolated vertices in vertex-deleted subgraphs. Amahashi and Kano \cite{AKo} proposed a characterization for
isolated $\frac{1}{k}$-tough graphs by virtue of $\{K_{1,j}:1\leq j\leq k\}$-factors. Zhou, Xu and Sun \cite{ZXS} proved some results on the existence of $\{K_{1,j}:1\leq j\leq k\}$-factors in graphs. Kano and Saito \cite{KS1}
verified a graph $G$ with $i(G-S)\leq\frac{1}{m}|S|$ for any subset $S\subseteq V(G)$ contains a $\{K_{1,j}:m\leq j\leq2m\}$-factor. Klopp and Steffen \cite{KS} provided some sufficient for graphs to possess special $\{K_{1,1},K_{1,2},\ldots,K_{1,k},C_m:m\geq3\}$-factors. Kano, Lu and Yu \cite{KLY} created a characterization for a graph with a $\{K_{1,1},K_{1,2},\ldots,K_{1,k},\mathcal{T}(2k+1)\}$-factor, where $\mathcal{T}(2k+1)$ is one
special family of trees. Wolf \cite{Wolf} characterized a graph with a $\{C_{2i+1},T:1\leq i<\frac{r}{k-r},T\in\mathcal{T}_{\frac{k}{r}}\}$-factor by using the number of isolated vertices, where $\mathcal{T}_{\frac{k}{r}}$
denotes the set of trees $T$ such that $i(T-S)\leq\frac{k}{r}|S|$ for all $S\subset V(T)$ and for every $e\in E(T)$ there exists a set $S^{*}\subset V(T)$ with $i((T-e)-S^{*})>\frac{k}{r}|S^{*}|$.

Recently, many scholars investigated the relationship between several spectral radii and component factors in graphs. O \cite{Os} posed an adjacency spectral radius condition for a connected graph with a $K_{1,1}$-factor. Zhou
and Zhang \cite{ZZ1} established a connection between a signless Laplacian spectral radius and a $K_{1,1}$-factor with a given property in a connected graph. Zhang and Lin \cite{ZLp} created a distance spectral radius condition
which ensures a connected graph and a connected bipartite graph to contain a $K_{1,1}$-factor, respectively. Li and Miao \cite{LM} proposed a sufficient condition to guarantee that a graph $G$ has a $P_{\geq2}$-factor by virtue
of its adjacency spectral radius. Zhou, Sun and Liu \cite{ZSL} put forward a distance signless Laplacian spectral radius condition for a connected graph to contain a $P_{\geq2}$-factor. Zhou, Zhang and Sun \cite{ZZS} obtained a
result on a $P_{\geq2}$-factor in a connected graph with respect to its $A_{\alpha}$-spectral radius. Wu \cite{Wc}, Liu, Fan and Shu \cite{LFS}, Zhou \cite{Zs1}, Wu, Zhou and Liu \cite{WZL} presented some adjacency spectral
conditions for the existence of tree-factors in connected graphs and connected bipartite graphs, respectively. Miao and Li \cite{ML} established a lower bound on the adjacency spectral radius for a connected graph $G$ to have a
$\{K_{1,j}:1\leq j\leq m\}$-factor and created an upper bound on the distance spectral radius for a connected graph $G$ to possess a $\{K_{1,j}:1\leq j\leq m\}$-factor. Wu \cite{Ws2} showed an adjacency spectral radius condition
for the existence of a $\{C_{2i+1},T:1\leq i<\frac{r}{k-r},T\in\mathcal{T}_{\frac{k}{r}}\}$-factor in a connected graph with minimum degree $\delta$, and claimed a signless Laplacian spectral radius condition for the existence
of a $\{C_{2i+1},T:1\leq i<\frac{r}{k-r},T\in\mathcal{T}_{\frac{k}{r}}\}$-factor in a connected graph with minimum degree $\delta$. Lv, Li and Xu \cite{LLX} gave two spectral conditions for a connected graph to contain a
$\{K_{1,1},C_{2i+1}:i\geq1\}$-factor. Zhou \cite{Zs} proposed some spectral conditions for a connected bipartite graph to have a special star-factor. Zhou and Liu \cite{ZL} provided two spectral conditions for the existence of
a $\{K_{1,j}:m\leq j\leq2m\}$-factor in a connected graph. We refer the reader to \cite{Ws1,LLY,ALY,Za,ZBW,ZZL} for other results on the connection between spectral radii and graph factors.

In this paper, we propose three sufficient conditions for a connected isolated tough graph to possess a $\{K_{1,j}:m\leq j\leq2m\}$-factor by means of the adjacency spectral radius, the signless Laplacian spectral radius and
the distance spectral radius.

\medskip

\noindent{\textbf{Theorem 1.1.}} Let $m$ and $b$ be two positive integers with $m\geq2$, and let $G$ be a connected isolated $\frac{mb-1}{b}$-tough graph of order $n$ with $n\geq\max\{(b^{2}+2b+1)m,(2b+2)m^{2}+(4b+1)m+2b-1\}$.
If
$$
\rho(G)\geq\rho(K_{mb-1}\vee(K_{n-(m+1)b+1}\cup bK_1)),
$$
then $G$ contains a $\{K_{1,j}:m\leq j\leq2m\}$-factor unless $G=K_{mb-1}\vee(K_{n-(m+1)b+1}\cup bK_1)$.

\medskip

\noindent{\textbf{Theorem 1.2.}} Let $m$ and $b$ be two positive integers with $m\geq2$, and let $G$ be a connected isolated $\frac{mb-1}{b}$-tough graph of order $n$ with $n\geq\max\{(b^{2}+2b+1)m+b^{2}-b,(2b+2)m^{2}+(4b+1)m+2b-1\}$.
If
$$
q(G)\geq q(K_{mb-1}\vee(K_{n-(m+1)b+1}\cup bK_1)),
$$
then $G$ contains a $\{K_{1,j}:m\leq j\leq2m\}$-factor unless $G=K_{mb-1}\vee(K_{n-(m+1)b+1}\cup bK_1)$.

\medskip

\noindent{\textbf{Theorem 1.3.}} Let $m$ and $b$ be two positive integers with $m\geq2$, and let $G$ be a connected isolated $\frac{mb-1}{b}$-tough graph of order $n$ with $n\geq\max\{2b^{2}m+b^{2}-b+1,3m^{2}b+3mb+m+4b\}$.
If
$$
\mu(G)\leq\mu(K_{mb-1}\vee(K_{n-(m+1)b+1}\cup bK_1)),
$$
then $G$ contains a $\{K_{1,j}:m\leq j\leq2m\}$-factor unless $G=K_{mb-1}\vee(K_{n-(m+1)b+1}\cup bK_1)$.

\medskip

\section{Some preliminaries}

Kano and Saito \cite{KS} proposed a sufficient condition for the existence of a $\{K_{1,j}:m\leq j\leq2m\}$-factor in a graph, where $m$ is an integer with $m\geq2$.

\medskip

\noindent{\textbf{Lemma 2.1}} (Kano and Saito \cite{KS}). Let $m$ be an integer with $m\geq2$, and let $G$ be a graph. If
$$
i(G-S)\leq\frac{1}{m}|S|
$$
for all subset $S\subseteq V(G)$, then $G$ has a $\{K_{1,j}:m\leq j\leq2m\}$-factor.

\medskip

\noindent{\textbf{Lemma 2.2}} (Li and Feng \cite{LF}). Let $H$ be a subgraph of a connected graph $G$. Then
$$
\rho(G)\geq\rho(H),
$$
with equality if and only if $G=H$.

\medskip

\noindent{\textbf{Lemma 2.3}} (Hong \cite{Ha}). Let $G$ be a graph of order $n$. Then
$$
\rho(G)\leq\sqrt{2e(G)-n+1},
$$
with equality occurring if and only if $G\in\{K_{1,n-1},K_n\}$.

\medskip

\noindent{\textbf{Lemma 2.4}} (Shen, You, Zhang and Li \cite{SYZL}). Let $H$ be a subgraph of a connected graph $G$. Then
$$
q(G)\geq q(H),
$$
where the equality occurs if and only if $G=H$.

\medskip

\noindent{\textbf{Lemma 2.5}} (Das \cite{Das}). Let $G$ be a graph with $n$ vertices. Then
$$
q(G)\leq\frac{2e(G)}{n-1}+n-2.
$$

\medskip

\noindent{\textbf{Lemma 2.6}} (Min\'c \cite{Mn}). Let $u$ and $v$ be two nonadjacent vertices of a connected graph $G$. Then
$$
\mu(G)>\mu(G+uv),
$$
where $G+uv$ denotes the graph obtained from $G$ by adding an edge to connect $u$ and $v$.

\medskip

By means of the Rayleigh quotient \cite{HJ}, the following lemma holds.

\medskip

\noindent{\textbf{Lemma 2.7.}} Let $G$ be a connected graph of order $n$. Then
$$
\mu(G)=\max_{X\neq\mathbf{0}}\frac{X^{T}\mathcal{D}(G)X}{X^{T}X}\geq\frac{\mathbf{1}^{T}\mathcal{D}(G)\mathbf{1}}{\mathbf{1}^{T}\mathbf{1}}=\frac{2W(G)}{n},
$$
where $\mathbf{1}=(1,1,\ldots,1)^{T}$.

\medskip

Let $M$ be an $n\times n$ real matrix, and let $\mathcal{N}=\{1,2,\ldots,n\}$. For a given partition $\pi:\mathcal{N}=\{\mathcal{N}_1,\mathcal{N}_2,\ldots,\mathcal{N}_r\}$ with
$\mathcal{N}=\mathcal{N}_1\cup\mathcal{N}_2\cup\ldots\cup\mathcal{N}_r$, the matrix $M$ can be correspondingly partitioned as
\begin{align*}
M=\left(
  \begin{array}{cccc}
    M_{11} & M_{12} & \cdots & M_{1r}\\
    M_{21} & M_{22} & \cdots & M_{2r}\\
    \vdots & \vdots & \ddots & \vdots\\
    M_{r1} & M_{r2} & \cdots & M_{rr}\\
  \end{array}
\right).
\end{align*}
The quotient matrix of the matrix $M$ based on the partition $\pi$ is defined by the matrix $M_{\pi}=(m_{ij})_{r\times r}$, where $m_{ij}$ denotes the average value of all row sum of the block $M_{ij}$. If for every pair $i,j$,
the block $M_{ij}$ has constant row sum $m_{ij}$, then the partition $\pi$ is called equitable.

\medskip

\noindent{\textbf{Lemma 2.8}} (You, Yang, So and Xi \cite{YYSX}). Let $M$ denote an $n\times n$ real matrix with an equitable partition $\pi$, and let $M_{\pi}$ be the corresponding quotient matrix. Then the eigenvalues
of $M_{\pi}$ are eigenvalues of $M$. Furthermore, if $M$ is nonnegative and irreducible, then the largest eigenvalues of $M$ and $M_{\pi}$ are equal.

\medskip

\section{The proof of Theorem 1.1}

\noindent{\it Proof of Theorem 1.1.} Suppose that a connected isolated $\frac{mb-1}{b}$-tough graph $G$ contains no $\{K_{1,j}:m\leq j\leq2m\}$-factor, where $m$ is an integer with $m\geq2$. Applying Lemma 2.1, there
exists some nonempty subset $S\subseteq V(G)$ satisfying
\begin{align}\label{eq:3.1}
|S|\leq m\cdot i(G-S)-1.
\end{align}
Let $|S|=s$ and $i(G-S)=i$. If $n\geq(m+1)i$, then $G$ is a spanning subgraph of $G_1=K_{mi-1}\vee(K_{n-(m+1)i+1}\cup iK_1)$. By means of Lemma 2.2, we obtain
\begin{align}\label{eq:3.2}
\rho(G)\leq\rho(G_1),
\end{align}
where the equality occurs if and only if $G=G_1$.

By virtue of the definition of isolated $\frac{mb-1}{b}$-tough graphs, we have
$$
\frac{mb-1}{b}\leq I(G)\leq\frac{|S|}{i(G-S)}=\frac{s}{i}.
$$
Combining this with \eqref{eq:3.1}, we possess
$$
bs\geq bmi-i\geq b(s+1)-i,
$$
which yields $i\geq b$.

If $i=b$, then we have $G_1=K_{mb-1}\vee(K_{n-(m+1)b+1}\cup bK_1)$. Together with \eqref{eq:3.2}, we conclude
$$
\rho(G)\leq\rho(K_{mb-1}\vee(K_{n-(m+1)b+1}\cup bK_1)),
$$
with equality if and only if $G=K_{mb-1}\vee(K_{n-(m+1)b+1}\cup bK_1)$, a contradiction. For $i\geq b+1$, using Lemma 2.3, we obtain
\begin{align}\label{eq:3.3}
\rho(G_1)\leq&\sqrt{2e(G_1)-n+1}\nonumber\\
=&\sqrt{2\binom{n-i}{2}+2i(mi-1)-n+1}\nonumber\\
=&\sqrt{(n-i)(n-i-1)+2i(mi-1)-n+1}\nonumber\\
=&\sqrt{(2m+1)i^{2}-(2n+1)i+n^{2}-2n+1}.
\end{align}
Let $f(i)=(2m+1)i^{2}-(2n+1)i+n^{2}-2n+1$. Since $n\geq(m+1)i$, we deduce $i\leq\frac{n}{m+1}$. Thus, $b+1\leq i\leq\frac{n}{m+1}$. A simple computation implies that
\begin{align*}
f(b+1)-f\left(\frac{n}{m+1}\right)=\frac{(n-mb-m-b-1)(n-(2b+2)m^{2}-(3b+2)m-b)}{(m+1)^{2}}\geq0,
\end{align*}
where the inequality occurs from the fact that $n\geq(2b+2)m^{2}+(4b+1)m+2b-1\geq(2b+2)m^{2}+(3b+2)m+b$. From the above discussion, we see
\begin{align}\label{eq:3.4}
f(b+1)\geq f(i)
\end{align}
for $b+1\leq i\leq\frac{n}{m+1}$. It follows from \eqref{eq:3.3}, \eqref{eq:3.4}, $m\geq2$, $b\geq1$ and $n\geq(b^{2}+2b+1)m$ that
\begin{align}\label{eq:3.5}
\rho(G_1)\leq&\sqrt{f(i)}\nonumber\\
\leq&\sqrt{f(b+1)}\nonumber\\
=&\sqrt{(2m+1)(b+1)^{2}-(2n+1)(b+1)+n^{2}-2n+1}\nonumber\\
=&\sqrt{(n-b-1)^{2}-2n+(2b^{2}+4b+2)m-b}\nonumber\\
\leq&\sqrt{(n-b-1)^{2}-2(b^{2}+2b+1)m+(2b^{2}+4b+2)m-b}\nonumber\\
<&n-b-1.
\end{align}

Notice that $K_{n-b}$ is a proper subgraph of $K_{mb-1}\vee(K_{n-(m+1)b+1}\cup bK_1)$. In terms of Lemma 2.2, we possess
\begin{align}\label{eq:3.6}
\rho(K_{mb-1}\vee(K_{n-(m+1)b+1}\cup bK_1))>\rho(K_{n-b})=n-b-1.
\end{align}
It follows from \eqref{eq:3.2}, \eqref{eq:3.5} and \eqref{eq:3.6} that
$$
\rho(G)\leq\rho(G_1)<n-b-1<\rho(K_{mb-1}\vee(K_{n-(m+1)b+1}\cup bK_1)),
$$
which contradicts $\rho(G)\geq\rho(K_{mb-1}\vee(K_{n-(m+1)b+1}\cup bK_1))$.

Let $n\leq(m+1)i-1$. Then $G$ is a spanning subgraph of $K_{n-i}\vee iK_1$ with $i\geq\lceil\frac{n+1}{m+1}\rceil$. Let $G_2=K_{n-\lceil\frac{n+1}{m+1}\rceil}\vee\lceil\frac{n+1}{m+1}\rceil K_1$. We easily see that
$K_{n-i}\vee iK_1$ is a spanning subgraph of $G_2$. According to Lemma 2.2, we possess
\begin{align}\label{eq:3.7}
\rho(G)\leq\rho(K_{n-i}\vee iK_1)\leq\rho(G_2),
\end{align}
where the equalities occur if and only if $G=G_2$. Notice that $2e(G_2)=(n-\lceil\frac{n+1}{m+1}\rceil)(n-\lceil\frac{n+1}{m+1}\rceil-1)+2\lceil\frac{n+1}{m+1}\rceil(n-\lceil\frac{n+1}{m+1}\rceil)
=(n-\lceil\frac{n+1}{m+1}\rceil)(n+\lceil\frac{n+1}{m+1}\rceil-1)<(n-\frac{n+1}{m+1})(n+\frac{n+1}{m+1})$. Combining this with Lemma 2.3, we get
\begin{align}\label{eq:3.8}
\rho(G_2)\leq&\sqrt{2e(G_2)-n+1}\nonumber\\
<&\sqrt{\Big(n-\frac{n+1}{m+1}\Big)\Big(n+\frac{n+1}{m+1}\Big)-n+1}.
\end{align}
Let $h(n)=(n-\frac{n+1}{m+1})(n+\frac{n+1}{m+1})-n+1$. Thus, we have
\begin{align*}
&(n-b-1)^{2}-h(n)\\
& \ =\frac{n^{2}-((2b+1)(m^{2}+2m+1)-2)n+(b^{2}+2b)(m^{2}+2m+1)+1}{(m+1)^{2}}\\
& \ >\frac{(b^{2}+2b)(m^{2}+2m+1)+1}{(m+1)^{2}}\\
& \ \ \ \ \ (\mbox{since} \ n\geq(2b+2)m^{2}+(4b+1)m+2b-1>(2b+1)m^{2}+(4b+2)m+2b-1)\\
& \ >0,
\end{align*}
which leads to
\begin{align}\label{eq:3.9}
\sqrt{h(n)}<n-b-1.
\end{align}
In terms of \eqref{eq:3.6}, \eqref{eq:3.7}, \eqref{eq:3.8} and \eqref{eq:3.9}, we conclude
$$
\rho(G)\leq\rho(G_2)<\sqrt{h(n)}<n-b-1<\rho(K_{mb-1}\vee(K_{n-(m+1)b+1}\cup bK_1)),
$$
which contradicts $\rho(G)\geq\rho(K_{mb-1}\vee(K_{n-(m+1)b+1}\cup bK_1))$. Theorem 1.1 is verified. \hfill $\Box$

\section{The proof of Theorem 1.2}

\noindent{\it Proof of Theorem 1.2.} Suppose that a connected isolated $\frac{mb-1}{b}$-tough graph $G$ contains no $\{K_{1,j}:m\leq j\leq2m\}$-factor, where $m\geq2$ is an integer. By Lemma 2.1, there exists some
nonempty subset $S\subseteq V(G)$ such that
\begin{align}\label{eq:4.1}
|S|\leq m\cdot i(G-S)-1.
\end{align}
Let $|S|=s$ and $i(G-S)=i$. If $n\geq(m+1)i$, then $G$ is a spanning subgraph of $G_1=K_{mi-1}\vee(K_{n-(m+1)i+1}\cup iK_1)$. According to Lemma 2.4, we conclude
\begin{align}\label{eq:4.2}
q(G)\leq q(G_1),
\end{align}
with equality occurring if and only if $G=G_1$.

In view of the definition of isolated $\frac{mb-1}{b}$-tough graphs, we have
\begin{align}\label{eq:4.3}
\frac{mb-1}{b}\leq I(G)\leq\frac{|S|}{i(G-S)}=\frac{s}{i}.
\end{align}
From \eqref{eq:4.1} and \eqref{eq:4.3}, we infer
$$
bs\geq bmi-i\geq b(s+1)-i,
$$
which leads to $i\geq b$.

If $i=b$, then $G_1=K_{mb-1}\vee(K_{n-(m+1)b+1}\cup bK_1)$. Using \eqref{eq:4.2}, we get
$$
q(G)\leq q(K_{mb-1}\vee(K_{n-(m+1)b+1}\cup bK_1)),
$$
with equality holding if and only if $G=K_{mb-1}\vee(K_{n-(m+1)b+1}\cup bK_1)$, a contradiction. If $i\geq b+1$, then it follows from Lemma 2.5 that
\begin{align}\label{eq:4.4}
q(G_1)\leq&\frac{2e(G_1)}{n-1}+n-2\nonumber\\
=&\frac{2\binom{n-i}{2}+2i(mi-1)}{n-1}+n-2\nonumber\\
=&\frac{(n-i)(n-i-1)+2i(mi-1)}{n-1}+n-2\nonumber\\
=&\frac{(2m+1)i^{2}-(2n+1)i+2n^{2}-4n+2}{n-1}.
\end{align}
Let $g(i)=(2m+1)i^{2}-(2n+1)i+2n^{2}-4n+2$. Recall that $i\geq b+1$ and $n\geq(m+1)i$. Then we deduce $b+1\leq i\leq\frac{n}{m+1}$. By a simple calculation, we obtain
\begin{align*}
g(b+1)-g\left(\frac{n}{m+1}\right)=\frac{(n-mb-m-b-1)(n-(2b+2)m^{2}-(3b+2)m-b)}{(m+1)^{2}}\geq0,
\end{align*}
where the inequality holds from the fact that $n\geq(2b+2)m^{2}+(4b+1)m+2b-1\geq(2b+2)m^{2}+(3b+2)m+b$, which implies
\begin{align}\label{eq:4.5}
g(b+1)\geq g(i)
\end{align}
for $b+1\leq i\leq\frac{n}{m+1}$. By means of \eqref{eq:4.4}, \eqref{eq:4.5}, $m\geq2$, $b\geq1$ and $n\geq(b^{2}+2b+1)m+b^{2}-b$, we deduce
\begin{align}\label{eq:4.6}
q(G_1)\leq&\frac{g(i)}{n-1}\nonumber\\
\leq&\frac{g(b+1)}{n-1}\nonumber\\
=&\frac{(2m+1)(b+1)^{2}-(2n+1)(b+1)+2n^{2}-4n+2}{n-1}\nonumber\\
=&2(n-b-1)-\frac{2n-2(b^{2}+2b+1)m-b^{2}+b}{n-1}\nonumber\\
\leq&2(n-b-1)-\frac{2((b^{2}+2b+1)m+b^{2}-b)-2(b^{2}+2b+1)m-b^{2}+b}{n-1}\nonumber\\
=&2(n-b-1)-\frac{b^{2}-b}{n-1}\nonumber\\
\leq&2(n-b-1).
\end{align}

Since $K_{n-b}$ is a proper subgraph of $K_{mb-1}\vee(K_{n-(m+1)b+1}\cup bK_1)$, it follows from Lemma 2.4 that
\begin{align}\label{eq:4.7}
q(K_{mb-1}\vee(K_{n-(m+1)b+1}\cup bK_1))>q(K_{n-b})=2(n-b-1).
\end{align}
By \eqref{eq:4.2}, \eqref{eq:4.6} and \eqref{eq:4.7}, we conclude
$$
q(G)\leq q(G_1)\leq2(n-b-1)<q(K_{mb-1}\vee(K_{n-(m+1)b+1}\cup bK_1)),
$$
which contradicts $q(G)\geq q(K_{mb-1}\vee(K_{n-(m+1)b+1}\cup bK_1))$.

Let $n\leq(m+1)i-1$. Then $G$ is a spanning subgraph of $K_{n-i}\vee iK_1$ with $i\geq\lceil\frac{n+1}{m+1}\rceil$. Let $G_2=K_{n-\lceil\frac{n+1}{m+1}\rceil}\vee\lceil\frac{n+1}{m+1}\rceil K_1$. Obviously,
$K_{n-i}\vee iK_1$ is a spanning subgraph of $G_2$. By virtue of Lemma 2.4, we infer
\begin{align}\label{eq:4.8}
q(G)\leq q(K_{n-i}\vee iK_1)\leq q(G_2),
\end{align}
with equalities if and only if $G=G_2$. Notice that $2e(G_2)=(n-\lceil\frac{n+1}{m+1}\rceil)(n-\lceil\frac{n+1}{m+1}\rceil-1)+2\lceil\frac{n+1}{m+1}\rceil(n-\lceil\frac{n+1}{m+1}\rceil)
=(n-\lceil\frac{n+1}{m+1}\rceil)(n+\lceil\frac{n+1}{m+1}\rceil-1)<(n-\frac{n+1}{m+1})(n+\frac{n+1}{m+1})$. Together with Lemma 2.5, we possess
\begin{align}\label{eq:4.9}
q(G_2)\leq&\frac{2e(G_2)}{n-1}+n-2\nonumber\\
<&\frac{(n-\frac{n+1}{m+1})(n+\frac{n+1}{m+1})}{n-1}+n-2\nonumber\\
=&2(n-b-1)-\frac{n^{2}-((2b+1)m^{2}+(4b+2)m+2b-1)n+2m^{2}b+4mb+2b+1}{(n-1)(m+1)^{2}}\nonumber\\
<&2(n-b-1)-\frac{2m^{2}b+4mb+2b+1}{(n-1)(m+1)^{2}}\nonumber\\
& \ \ \ (\mbox{since} \ n\geq(2b+2)m^{2}+(4b+1)m+2b-1>(2b+1)m^{2}+(4b+2)m+2b-1)\nonumber\\
<&2(n-b-1).
\end{align}
Using \eqref{eq:4.7}, \eqref{eq:4.8} and \eqref{eq:4.9}, we obtain
$$
q(G)\leq q(G_2)<2(n-b-1)<q(K_{mb-1}\vee(K_{n-(m+1)b+1}\cup bK_1)),
$$
which contradicts $q(G)\geq q(K_{mb-1}\vee(K_{n-(m+1)b+1}\cup bK_1))$. Theorem 1.2 is proved. \hfill $\Box$

\section{The proof of Theorem 1.3}

\noindent{\it Proof of Theorem 1.3.} Suppose that a connected isolated $\frac{mb-1}{b}$-tough graph $G$ contains no $\{K_{1,j}:m\leq j\leq2m\}$-factor, where $m\geq2$ is an integer. According to Lemma 2.1, we have
\begin{align}\label{eq:5.1}
|S|\leq m\cdot i(G-S)-1
\end{align}
for some nonempty subset $S\subseteq V(G)$. Let $|S|=s$ and $i(G-S)=i$. The following proof will be divided into two cases in terms of the value of $n$.

\noindent{\bf Case 1.} $n\geq(m+1)i$.

Clearly, $G$ is a spanning subgraph of $G_1=K_{mi-1}\vee(K_{n-(m+1)i+1}\cup iK_1)$. Using Lemma 2.6, we conclude
\begin{align}\label{eq:5.2}
\mu(G)\geq\mu(G_1),
\end{align}
where the equality follows if and only if $G=G_1$.

From the definition of isolated $\frac{mb-1}{b}$-tough graphs, we get
$$
\frac{mb-1}{b}\leq I(G)\leq\frac{|S|}{i(G-S)}=\frac{s}{i}.
$$
Together with \eqref{eq:5.1}, we deduce
$$
bs\geq bmi-i\geq b(s+1)-i,
$$
which leads to $i\geq b$.

If $i=b$, then we conclude $G_1=K_{mb-1}\vee(K_{n-(m+1)b+1}\cup bK_1)$. By \eqref{eq:5.2}, we obtain
$$
\mu(G)\geq\mu(K_{mb-1}\vee(K_{n-(m+1)b+1}\cup bK_1)),
$$
where the equality holds if and only if $G=K_{mb-1}\vee(K_{n-(m+1)b+1}\cup bK_1)$, a contradiction. In what follows, we shall consider $i\geq b+1$.

Recall that $G_1=K_{mi-1}\vee(K_{n-(m+1)i+1}\cup iK_1)$. The quotient matrix of $\mathcal{D}(G_1)$ by virtue of the partition $V(G_1)=V(K_{mi-1})\cup V(K_{n-(m+1)i+1})\cup V(iK_1)$ can be given by
\begin{align*}
B_1=\left(
  \begin{array}{ccc}
    mi-2 & n-(m+1)i+1 & i\\
    mi-1 & n-(m+1)i & 2i\\
    mi-1 & 2n-(2m+2)i+2 & 2i-2\\
  \end{array}
\right),
\end{align*}
and so its characteristic polynomial is
\begin{align*}
\varphi_{B_1}(x)=&x^{3}-(i+n-4)x^{2}+((3m+2)i^{2}-(2n+4)i-3n+5)x\\
&-m(m+1)i^{3}+(mn+5m+3)i^{2}-(3n+4)i-2n+2.
\end{align*}
By means of Lemma 2.8 and the equitable partition $V(G_1)=V(K_{mi-1})\cup V(K_{n-(m+1)i+1})\cup V(iK_1)$, the largest root of $\varphi_{B_1}(x)=0$ equals $\mu(G_1)$.

Write $G_*=K_{mb-1}\vee(K_{n-(m+1)b+1}\cup bK_1)$. The quotient matrix of $\mathcal{D}(G_*)$ with respect to the partition $V(G_*)=V(K_{mb-1})\cup V(K_{n-(m+1)b+1})\cup V(bK_1)$ equals
\begin{align*}
B_*=\left(
  \begin{array}{ccc}
    mb-2 & n-(m+1)b+1 & b\\
    mb-1 & n-(m+1)b & 2b\\
    mb-1 & 2n-(2m+2)b+2 & 2b-2\\
  \end{array}
\right).
\end{align*}
The characteristic polynomial of $B_*$ is
\begin{align*}
\varphi_{B_*}(x)=&x^{3}-(b+n-4)x^{2}+((3m+2)b^{2}-(2n+4)b-3n+5)x\\
&-m(m+1)b^{3}+(mn+5m+3)b^{2}-(3n+4)b-2n+2.
\end{align*}
In view of Lemma 2.8 and the equitable partition $V(G_*)=V(K_{mb-1})\cup V(K_{n-(m+1)b+1})\cup V(bK_1)$, the largest root, say $\mu_*$, of $\varphi_{B_*}(x)=0$ equals $\mu(G_*)$. That is to say, $\mu(G_*)=\mu_*$
and $\varphi_{B_*}(\mu_*)=0$.

In terms of Lemma 2.7 and $n\geq2b^{2}m+b^{2}-b+1$, we possess
\begin{align}\label{eq:5.3}
\mu_*=&\mu(G_*)\nonumber\\
\geq&\frac{2W(G_*)}{n}\nonumber\\
=&\frac{n^{2}+(2b-1)n-2b^{2}m-b^{2}+b}{n}\nonumber\\
=&\frac{n^{2}+(2b-2)n+n-2b^{2}m-b^{2}+b}{n}\nonumber\\
>&n+2b-2.
\end{align}

Recall that $\varphi_{B_*}(\mu_*)=0$. By plugging the value $\mu_*$ into $x$ of $\varphi_{B_1}(x)-\varphi_{B_*}(x)$, we possess
\begin{align}\label{eq:5.4}
\varphi_{B_1}(\mu_*)=\varphi_{B_1}(\mu_*)-\varphi_{B_*}(\mu_*)=(i-b)\beta(\mu_*),
\end{align}
where $\beta(\mu_*)=-\mu_*^{2}+((3m+2)i+(3m+2)b-2n-4)\mu_*-m(m+1)i^{2}+(mn-bm^{2}-bm+5m+3)i-m(m+1)b^{2}+(mn+5m+3)b-3n-4$. Notice that
$$
\frac{(3m+2)i+(3m+2)b-2n-4}{2}<n+2b-2<\mu_*
$$
by \eqref{eq:5.3}, $i\geq b+1$, $n\geq(m+1)i$ and $n\geq3m^{2}b+3mb+m+4b$. Hence, we get
\begin{align*}
\beta(\mu_*)<&\beta(n+2b-2)\\
=&-(n+2b-2)^{2}+((3m+2)i+(3m+2)b-2n-4)(n+2b-2)-m(m+1)i^{2}\\
&+(mn-bm^{2}-bm+5m+3)i-m(m+1)b^{2}+(mn+5m+3)b-3n-4\\
=&-m(m+1)i^{2}+(4mn+2n-bm^{2}+5bm-m+4b-1)i\\
&-3n^{2}+(4mb-6b+1)n-m^{2}b^{2}+5mb^{2}-mb-b\\
\triangleq&\gamma(i).
\end{align*}
Notice that
$$
\frac{4mn+2n-bm^{2}+5bm-m+4b-1}{2m(m+1)}>\frac{n}{m+1}\geq i\geq b+1.
$$
Then we have
\begin{align}\label{eq:5.5}
\beta(\mu_*)<&\gamma(i)\nonumber\\
\leq&\gamma\left(\frac{n}{m+1}\right)\nonumber\\
=&-m(m+1)\left(\frac{n}{m+1}\right)^{2}+(4mn+2n-bm^{2}+5bm-m+4b-1)\left(\frac{n}{m+1}\right)\nonumber\\
&-3n^{2}+(4mb-6b+1)n-m^{2}b^{2}+5mb^{2}-mb-b\nonumber\\
=&\frac{1}{m+1}\Big(-n^{2}+(3m^{2}b+3mb-2b)n-m^{3}b^{2}+4m^{2}b^{2}-m^{2}b+5mb^{2}-2mb-b\Big)\nonumber\\
<&\frac{1}{m+1}\Big(-(3m^{2}b+3mb-b)^{2}+(3m^{2}b+3mb-2b)(3m^{2}b+3mb-b)-m^{3}b^{2}\nonumber\\
&+4m^{2}b^{2}-m^{2}b+5mb^{2}-2mb-b\Big) \ \ \ \ \ \ (\mbox{since} \ n>3m^{2}b+3mb-b)\nonumber\\
=&\frac{1}{m+1}\Big(-(m^{3}-m^{2}-2m+1)b^{2}-(m^{2}+2m-1)b\Big)\nonumber\\
<&0 \ \ \ \ \ \ \ (\mbox{since} \ b\geq1 \ \mbox{and} \ m\geq2).
\end{align}
Using \eqref{eq:5.4}, \eqref{eq:5.5} and $i\geq b+1$, we deduce
$$
\varphi_{B_1}(\mu_*)=(i-b)\beta(\mu_*)<0,
$$
which leads to
\begin{align}\label{eq:5.6}
\mu(G_1)>\mu_*=\mu(G_*)=\mu(K_{mb-1}\vee(K_{n-(m+1)b+1}\cup bK_1).
\end{align}
According to \eqref{eq:5.2} and \eqref{eq:5.6}, we conclude
$$
\mu(G)\geq\mu(G_1)>\mu(K_{mb-1}\vee(K_{n-(m+1)b+1}\cup bK_1),
$$
which contradicts $\mu(G)\leq\mu(K_{mb-1}\vee(K_{n-(m+1)b+1}\cup bK_1)$.

\noindent{\bf Case 2.} $n\leq(m+1)i-1$.

Obviously, $G$ is a spanning subgraph of $K_{n-i}\vee iK_1$ with $i\geq\lceil\frac{n+1}{m+1}\rceil$. Let $G_2=K_{n-\lceil\frac{n+1}{m+1}\rceil}\vee\lceil\frac{n+1}{m+1}\rceil K_1$. We easily see that
$K_{n-i}\vee iK_1$ is a spanning subgraph of $G_2$. Applying Lemma 2.6, we have
\begin{align}\label{eq:5.7}
\mu(G)\geq\mu(K_{n-i}\vee iK_1)\geq\mu(G_2),
\end{align}
with equalities if and only if $G=G_2$. By the partition $V(G_2)=V(K_{n-\lceil\frac{n+1}{m+1}\rceil})\cup V(\lceil\frac{n+1}{m+1}\rceil K_1)$, the quotient matrix of $\mathcal{D}(G_2)$ is
\begin{align*}
B_2=\left(
  \begin{array}{ccc}
    n-\lceil\frac{n+1}{m+1}\rceil-1 & \lceil\frac{n+1}{m+1}\rceil\\
    n-\lceil\frac{n+1}{m+1}\rceil & 2\lceil\frac{n+1}{m+1}\rceil-2\\
  \end{array}
\right).
\end{align*}
By means of a simple calculation, we know that the characteristic polynomial of $B_2$ equals
$$
\varphi_{B_2}(x)=x^{2}-\Big(n+\Big\lceil\frac{n+1}{m+1}\Big\rceil-3\Big)x+n\Big\lceil\frac{n+1}{m+1}\Big\rceil-\Big\lceil\frac{n+1}{m+1}\Big\rceil^{2}-2n+2.
$$
In view of Lemma 2.8 and the equitable partition $V(G_2)=V(K_{n-\lceil\frac{n+1}{m+1}\rceil})\cup V(\lceil\frac{n+1}{m+1}\rceil K_1)$, $\mu(G_2)$ is the largest root of $\varphi_{B_2}(x)=0$.

Recall that $G_*=K_{mb-1}\vee(K_{n-(m+1)b+1}\cup bK_1)$, $\varphi_{B_*}(x)=x^{3}-(b+n-4)x^{2}+((3m+2)b^{2}-(2n+4)b-3n+5)x-m(m+1)b^{3}+(mn+5m+3)b^{2}-(3n+4)b-2n+2$ and $\mu_*=\mu(G_*)$ is the largest root of
$\varphi_{B_*}(x)=0$. Notice that
$$
\varphi_{B_*}'(x)=3x^{2}-2(b+n-4)x+(3m+2)b^{2}-(2n+4)b-3n+5
$$
and $\frac{b+n-4}{3}<n-b-2$. For $x\geq n-b-2$, we have
\begin{align*}
\varphi_{B_*}'(x)>&\varphi_{B_*}'(n-b-2)\\
=&n^{2}-(8b+3)n+(3m+7)b^{2}+4b+1\\
>&0 \ \ \ \ \ (\mbox{since} \ m\geq2, \ b\geq1 \ \mbox{and} \ n\geq3m^{2}b+3mb+m+4b),
\end{align*}
which implies that $\varphi_{B_*}(x)$ is increasing for $x\geq n-b-2$. By a simple computation, we obtain
$$
\varphi_{B_*}(n+3b-1)=(4m+8)b^{2}n+(-m^{2}+8m+24)b^{3}+(2m+4)b^{2}-b>0,
$$
which leads to
\begin{align}\label{eq:5.8}
\mu(G_*)<n+3b-1.
\end{align}
By a direct calculation, we get
\begin{align*}
\varphi_{B_2}(n+3b-1)=&3bn-(3b-1)\Big\lceil\frac{n+1}{m+1}\Big\rceil-\Big\lceil\frac{n+1}{m+1}\Big\rceil^{2}+9b^{2}+3b\\
\leq&3bn-(3b-1)\Big(\frac{n+1}{m+1}\Big)-\Big(\frac{n+1}{m+1}\Big)^{2}+9b^{2}+3b\\
=&\frac{-n^{2}+(3bm^{2}+3bm+m-1)n+(9b^{2}+3b)m^{2}+(18b^{2}+3b+1)m+9b^{2}+1}{(m+1)^{2}}\\
\leq&\frac{-3m^{2}b^{2}+6mb^{2}-4mb-7b^{2}-4b+1}{(m+1)^{2}}\\
& \ \ \ \ \ (\mbox{since} \ n\geq3m^{2}b+3mb+m+4b)\\
<&0 \ \ \ \ \ (\mbox{since} \ m\geq2 \ \mbox{and} \ b\geq1),
\end{align*}
which leads to
\begin{align}\label{eq:5.9}
\mu(G_2)>n+3b-1.
\end{align}
It follows from \eqref{eq:5.7}, \eqref{eq:5.8} and \eqref{eq:5.9} that
$$
\mu(G)\geq\mu(G_2)>n+3b-1>\mu(G_*)=\mu(K_{mb-1}\vee(K_{n-(m+1)b+1}\cup bK_1)),
$$
which contradicts $\mu(G)\leq\mu(K_{mb-1}\vee(K_{n-(m+1)b+1}\cup bK_1))$. Theorem 1.3 is proved. \hfill $\Box$

\section*{Data availability statement}

My manuscript has no associated data.

\section*{Declaration of competing interest}

The authors declare that they have no conflicts of interest to this work.

\section*{Acknowledgments}

This work was supported by the Natural Science Foundation of Jiangsu Province (Grant No. BK20241949). Project ZR2023MA078 supported by Shandong Provincial Natural Science Foundation.

\end{document}